\documentclass[12pt]{article}
\usepackage{amssymb}
\def\ad{\mbox{ad}}
\def\qed{\ \ \ifhmode\unskip\nobreak\fi\ifmmode\ifinner
         \else\hskip5pt\fi\fi
 \hbox{\hskip5pt\vrule width4pt height6pt depth1.5pt\hskip 1 pt}}
\def\a{\alpha}

\def\d{\delta}

\def\l{\lambda}

\def\Vir{\mbox{Vir}}

\def\si{\sigma}
\def\sc{\scriptstyle}
\def\ssc{\scriptscriptstyle}
\def\Z{\mathbb{Z}}\def\C{\mathbb{C}}
\def\dis{\displaystyle}
\def\cl{\centerline}

\def\nl{\newline}
\def\ol{\overline}

\def\wh{\widehat}

\def\bs{\backslash}

\def\rb{\raisebox}
\def\vs{\vspace*}
\def\ra{\rangle}
\def\la{\langle}
\def\ni{\noindent}
\def\hi{\hangindent}
\def\ha{\hangafter}
\def\WW{{\cal W}}

\def\HH{{\cal H}}
\def\es{\varepsilon}
\def\DD{{\mathcal{D}}}
\def\LIE{\rb{-1.pt}{\mbox{$\wh{\DD^N}$}}}
\def\LI{\rb{-1.pt}{\mbox{$\wh{\DD^1}$}}}
\def\cc{{ C}}
\begin{document}
\cl {\large\bf Classification of Quasifinite Modules over
        Lie Algebras}\cl{\large\bf
of Matrix Differential Operators on the Circle} \vs{10pt}\par
\cl{(appeared in {\it Proc.~Amer.~Math.~Soc.}, 
{\bf133} (2005), 1949-1957)}
 \vs{10pt}\par
\cl{{\bf Yucai Su}\footnote{Supported by a NSF grant 10171064 of
China and 
two grants ``Excellent Young Teacher Program'' and
``Trans-Century Training Programme Foundation for the Talents''
from Ministry of Education of
China.\vs{4pt}\nl\hspace*{4.5ex}{\it Mathematics Subject
Classification (1991):} 17B10, 17B65, 17B66, 17B68}} \cl{\small\it
Department of Mathematics, Shanghai Jiaotong University, Shanghai
200030, P.~R.~China {\rm and}}\cl{\small\it Department of
Mathematics, Harvard University, Cambridge, MA 02138,
USA}\cl{\small\it Email: ycsu@sjtu.edu.cn} \vs{10pt}\par {\small
{\bf Abstract} \ We prove that an irreducible quasifinite module
over the central extension of the Lie algebra of $N\times
N$-matrix differential operators on the circle is either a highest
or lowest weight module or else a module of the intermediate
series. Furthermore, we give a complete classification of
indecomposable uniformly
bounded modules.}%
\vs{10pt}\par \cl{1. \ INTRODUCTION}
\par
Kac [5] introduced the notion of conformal algebras. Conformal
algebras play important roles in quantum field theory and vertex
operator algebras (e.g.~[5]), whose study has drawn much attention
in the literature (e.g.~[1, 5-7, 16-19]). It is pointed out in [7]
that the infinite rank associative conformal algebra ${\rm
Cend}_N$ and Lie conformal algebra $gc_N$ (the {\it general Lie
conformal algebra}) play the same important roles in the theory of
conformal algebras as ${\rm End}_N$ and $gl_N$ play in the theory
of associative and Lie algebras.
\par
There is a one to one correspondence between Lie conformal
algebras and maximal formal distribution Lie
algebras [1, 6, 7]. The Lie algebra $\DD^N$
of $N\times N$-matrix differential operators on the circle is a
formal distribution Lie algebra associated to the general Lie
conformal algebra $gc_N$. Let
\LIE\ be the universal central extension of
$\DD^N$. In particular when $N=1$, \LI\ is also known as the $W$-infinity
algebra $\WW_{1+\infty}$.
Thus one may expect that the representation
theory of $\DD^N$ and \LIE\ is very
interesting and important (e.g.~[2, 4, 8-10, 13]). %
\par%
As is pointed in [8, 10], although \LIE\ is a $\Z$-graded Lie
algebra, each grading space is still infinite dimensional
(cf.~(2.4)), the classification of quasifinite modules is thus a
nontrivial problem. The classification of irreducible quasifinite
highest weight modules over \LIE\ was given by Kac and Radul [8]
for case $N=1$. For the general $N$, the classification was obtained
by Boyallian, Kac, Liberati and Yan [2]. 
In [13], the author obtained the classification of the
irreducible quasifinite modules and indecomposable uniformly
bounded modules over \LI. We would like to take this chance to
point out that there is a slight gap in the proof of Proposition
2.2 of [13], this gap has been filled in this paper, see Remark
3.4. 

In this paper, we generalize the results in [13] to the
general case: Precisely, we obtain that an irreducible quasifinite
module over \LIE\ (thus also over $\DD^N$) is either a highest or
lowest weight module or else a module of the intermediate series.
Furthermore, we give a complete classification of indecomposable
uniformly bounded modules (Theorem 2.2). %
\vskip 10pt\par%
\cl{2. \ NOTATION AND MAIN THEOREM }%
\par
Let $N\ge1$ be an integer. Let $\C[t,t^{-1}]$ be the Laurent
polynomial algebra over the variable $t$, let $\C[D]$ be the
polynomial algebra over $D=t{d\over dt}$, and let $gl_N$ be the
space of $N\times N$ matrices. The {\it Lie algebra $\DD^N$ of
$N\times N$-matrix differential operators on the circle} is the
tensor product space $\DD^N=\C[t,t^{-1}]\otimes\C[D]\otimes gl_N$,
spanned by $\{t^i D^jA\,|\ i\in\Z,\,j\in\Z_+,\,A\in gl_N\}$, with
the Lie bracket:
$$
[t^i D^jA,t^k D^lB]= (t^i D^jA)\cdot(t^k D^lB) -(t^k
D^lB)\cdot(t^i D^jA),%
\eqno(2.1)$$%
\vs{-4pt}and
$$ (t^iD^jA)\cdot(t^kD^lB)
=t^{i+k}(D+k)^jD^lAB=\sum_{s=0}^j\left(\!\begin{array}{c}j\\
s\end{array}\!\right) k^s t^{j+k}D^{j+l-s}AB,
\vs{-4pt}\eqno(2.2)$$%
for $i,j\in\Z,\,j,l\in\Z_+$, where
$(^i_s)={i(i-1)\cdots(i-s+1)\over s!}$ if $s\ge0$ and $(^i_s)=0$
if $s<0$, is the binomial coefficient. The associative algebra
with the underlined vector space $\DD^N$ and the product (2.2) is
denoted by $\DD^N_{\rm as}$.
\par
It is proved in [11] that $\DD^N$ has a unique nontrivial central
extension. The {\it universal central extension \LIE\ of $\DD^N$}
is defined as follows (cf.~[2]): The Lie bracket (2.1) is replaced
by (cf.\vs{-4pt}~(3.12)) $$
\begin{array}{ll}
[t^i[D]_jA,t^k[D]_lB]=\!\!\!\!& (t^i[D]_jA)\cdot(t^k[D]_lB)
-(t^k[D]_lB)\cdot(t^i[D]_jA)\vs{4pt}\\
&+\,\d_{i,-k}(-1)^jj!l!
\biggl(\!\begin{array}{c}i+j\\
j+l+1\end{array}\!\biggr){\rm tr}(AB)\cc, \end{array}%
 \vs{-7pt}\eqno(2.3)$$%
for $i,j\in\Z,\,j,l\in\Z_+$, where $[D]_j=D(D-1)
\cdots(D-j+1)=t^j({d\over dt})^j$, and ${\rm tr}(A)$ is the trace
of the matrix $A$, and where $\cc$ is a nonzero central element of
\LIE.
\par%
For $m,n\in\Z$, we denote $[m,n]=\{m,m+1,...,n\}$. Let
$\{E_{p,q}\,|\,p,q\in[1,N]\}$ be the standard basis of $gl_N$,
where $E_{p,q}$ is the matrix with entry $1$ at $(p,q)$ and $0$
otherwise. Then \LIE\ has the {\it principal $\Z$-gradation}
$\LIE= \oplus_{i\in\Z}(\LIE)_i$ with the grading space (cf.~[2])
$$
(\LIE)_i={\rm span}\{t^k D^j
E_{p,q}\,|\,k\in\Z,\,j\in\Z_+,\,p,q\in[1,N],\,
kN+p-q=i\}\oplus\d_{i,0}\C{\ssc\,}\cc,
\eqno(2.4)$$%
for $i\in\Z$. In particular, $(\LIE)_0={\rm
span}\{D^jE_{p,p},\cc\,|\,j\in\Z_+,p\in[1,N]\}$ is a commutative
subalgebra. Note that $t$ has degree $N$.
\par
{D$\sc\rm EFINITION$ 2.1.} \ A \LIE-module (or a $\DD^N_{\rm
as}$-module) $V$ is called a {\it quasifinite module} (e.g.~[2])
if $V=\oplus_{j\in\Z}V_j$ is $\Z$-graded such that
$(\LIE)_iV_j\subset V_{i+j}$ and ${\rm dim\,}V_j<\infty$ for
$i,j\in\Z$. This is equivalent to saying that a quasifinite module
is a module having finite dimensional generalized weight spaces
with respect to the commutative subalgebra $(\LIE)_0$. A
quasifinite module $V$ is called a {\it module of the intermediate
series} if ${\rm dim\ssc\,}V_j\le1$ for $j\in\Z$. It is called
 a {\it uniformly
bounded module} if there exists an integer $K>0$ such that ${\rm
dim\ssc\,}V_j\le K$ for $j\in\Z$. A module $V$ is
a {\it trivial module} if $\LIE$ acts trivially on $V$
(i.e., $\LIE \,V=0$).
\par
Clearly a $\DD^N_{\rm as}$-module is also a $\DD^N$-module (but
the converse is not necessarily true), and a \mbox{$\DD^N$-module} is a
\LIE-module (with central charge 0). Thus it suffices to consider
\LIE-modules.
\par%
We shall define 2 families of modules $V(\a),\ol V(\a)$,
$\a\in\C$, of the intermediate series \mbox{below}. For a fixed
$\a\in\C$, the obvious representation of \LIE\ (with trivial
action of the central element $\cc$) on the space
$V(\a)=\C^N[t,t^{-1}]t^\a$ defines an irreducible module $V(\a)$.
Let $\{\es_p=(\d_{p1},...,\d_{pN})^{\rm T}\,|\,p\in[1,N]\}$ be the
standard basis of $\C^N$, where the superscript ``T'' means the
transpose of vectors or matrices (thus we regard elements of
$\C^N$ as column vectors). Then action of \LIE\ on $V(\a)$
\vs{-2pt}is
$$
(t^iD^jE_{p,q})(t^{k+\a}\es_r)=\d_{q,r}(k+\a)^jt^{i+k+\a}\es_p, %
\vs{-2pt}\eqno(2.5)$$%
for $i,k\in\Z,\,j\in\Z_+,\,p,q,r\in[1,N]$. For $j\in\Z$,
\vs{-2pt}let
$$
V(\a)_j=\C t^{k+\a}\es_r, %
\vs{-2pt}\eqno(2.6)$$%
where $k\in\Z,\,r\in[1,N]$ are unique such that $j+1=kN+r,$ then
$V(\a)=\oplus_{j\in\Z}V(\a)_j$ is a $Z$-graded space such that
${\rm dim\sc\,}V(\a)_j=1$ for $j\in\Z$. Thus $V(\a)$ is a module
of the intermediate series.
\par
For $v\in\C^N,\,k\in\Z$, denote $v_k=t^{k+\a}v\in V(\a)$ (note
that $v_k$ is in general not a homogeneous vector). For $A\in
gl_N$, define $Av_k=(Av)_k$, where $Av$ is the natural action of
$A$ on $v$ defined linearly by $E_{p,q}\es_r=\d_{q,r}\es_p$ (i.e.,
the action is defined by the matrix-vector multiplication). Then
(2.5) can be rewritten \vs{-2pt}as
$$
(t^iD^jA) v_k=(k+\a)^jAv_{i+k},
\vs{-2pt}\eqno(2.7)$$%
for $i,k\in\Z,\,j\in\Z_+,\,A\in gl_N,\,v\in\C^N$. Clearly $V(\a)$
is also a $\DD^N_{\rm as}$-module.
\par
 Note that (cf.~[15]) there exists a Lie algebra
isomorphism $\si:\DD^N\cong\DD^N$ such \vs{-2pt}that
$$
\si(t^i D^jA)=(-1)^{j+1}t^i(D+i)^jA^{\rm T},%
\vs{-2pt}\eqno(2.8)$$%
for $i\in\Z,\,j\in\Z_+,\,A\in gl_N$. Using this isomorphism, we
have another \LIE-module $\ol V(\a)$ (again with trivial action of
$\cc$), called the {\it twisted module of $V(\a)$}, defined
\vs{-2pt}by
$$
(t^iD^jA)  v_k=(-1)^{j+1}(i+k+\a)^jA^{\rm T}  v_{i+k},
\vs{-4pt}\eqno(2.9)$$%
for $i,k\in\Z,\,j\in\Z_+,\,A\in gl_N,\,  v\in\C^N$. In fact, $\ol
V(\a)$ is the dual module of $V(-\a)$: If we define a bilinear
form on $\ol V(\a)\times V(-\a)$ by $\la
t^{i+\a}\es_p,t^{j-\a}\es_q\ra=\d_{i+j,0}\d_{p,q}$, then
$$
\la x  v,v\ra=-\la  v,xv\ra,%
$$
for $x\in\LIE,\,  v\in\ol V(\a),\,v\in V(-\a)$. Note that a
$\Z$-gradation of $\ol V(\a)$ can be defined by (2.6) with $k,r$
satisfying the relation $j=kN+N-r$. Obviously, $\ol V(\a)$ is not
a $\DD^N_{\rm as}$-module.
\par%
Now we can generalize the above modules  $V(\a)$ and $\ol V(\a)$
as follows: Let $m>0$ be an integer, let $\a$ be an indecomposable
linear transformation on $\C^m$ (thus up to equivalences, $\a$ is
uniquely determined by its unique eigenvalue $\l$). Let
$gl_N$, $\a$ act on $\C^N\otimes\C^m$ defined by $A(u\otimes
v)=Au\otimes v$, $\a(u\otimes v)=u\otimes \a v$ for $A\in
gl_N,\,u\in\C^N,\,u\in\C^m$. Then by allowing $v$ to be
in $\C^N\otimes\C^m$ in (2.7) and (2.9), we obtain 2 families of
indecomposable uniformly bounded modules
$V(m,\a)=V(\a)\otimes\C^m$, $\ol V(m,\a)=\ol V(\a)\otimes\C^m$.
\par
The main result of this paper is the following theorem.
\par {T$\sc\rm HEOREM$ 2.2.} \ {\it (i) An irreducible quasifinite module over \LIE\ is
either a highest or lowest weight module or else a module of the
intermediate series. A nontrivial module of the intermediate
series is a module $V(\a)$ or $\ol V(\a)$ for some $\a\in\C$. %
\vs{-4pt}\par%
(ii) A nontrivial indecomposable uniformly bounded module over
\LIE\ is a module $V(m,\a)$ or $\ol V(m,\a)$ for some
$m\in\Z_+\bs\{0\}$ and some indecomposable linear transformation
$\a$ of $\C^m$.}
\par Thus in particular we obtain that a nontrivial indecomposable uniformly
bounded module over $\DD^N$ is simply a $\DD^N_{\rm as\,}$-module
(if the central element $I=t^0D^0I$ acts by $1$, where $I$ is the
$N\times N$ identity matrix) or its twist (if $I$ acts by $-1$),
and that there is an equivalence between the category of uniformly
bounded $\DD^N_{\rm as}$-modules without the trivial composition
factor and the category of linear transformations on finite
dimensional vector spaces.
\par
Since irreducible quasifinite highest weight modules have been
classified in [2] and irreducible lowest weight modules are simply
dual modules of irreducible highest weight modules, this theorem
in fact classifies all irreducible quasifinite modules over \LIE\
and over
$\DD^N$. %
 \par%
Note that in the language of conformal algebras, this theorem
in particular also gives proofs of Theorems 6.1 and 6.2 of [7] on
the classification of finite indecomposable modules over the
conformal algebras ${\rm Cend}_N$ and $gc_N$.
\par%
 The analogous results of this theorem for affine Lie
algebras, the Virasoro algebra, higher rank Virasoro algebras and
Lie algebras of Weyl type were obtained in [3, 12-14].%
\par
Note that the space $\HH={\rm
span}\{\cc,D,E_{p,p}\,|\,p\in[1,N]\}$ is a Cartan subalgebra of
\LIE\ and that the definition of quasifiniteness does not require
that $V$ is a {\it weight module} (i.e., the actions of elements
of $\HH$ on $V$ are diagonalizable). If we require $V$ to be a
weight module, then the linear transformation $\a$ is
diagonalizable, and thus all uniformly bounded modules are
completely reducible.
\par
We shall prove the above theorem in the next section.%
 \vskip 10pt\par
\cl{3. \ PROOF OF THEOREM 2.2}
\vs{0pt}\par%
We shall keep notations of the previous section.
We denote by $I$ the $N\times N$ identity matrix. When the context
is clear, we often omit the symbol $I$; for instance,
$t=t^1D^0I\in\DD^N$. Denote $\HH={\rm span}
\{\cc,D,E_{p,p}\,|\,p\in[1,N]\}$, a Cartan subalgebra of \LIE.
\par%
For any $\Z$-graded vector space $U$, we use notations $U_+,\,U_-$
and $U_{[p,q)}$ to denote the subspaces spanned by 
elements of degree $k$ with $k>0,\,k<0$ and $p\le k< q$
respectively. Then \LIE\ has a triangular decomposition
$\LIE=(\LIE)_+\oplus(\LIE)_0\oplus(\LIE)_-$. Observe that
$(\LIE)_+$ is generated by the adjoint action of $t$ on
$(\LIE)_{[0,N)}$ and that $\ad_t$ is locally nilpotent on \LIE.
\par
Suppose $V=\oplus_{i\in\Z}V_i$ is a quasifinite module over $\WW$.
A homogeneous nonzero vector $v\in V$ is called a {\it highest}
(resp., {\it lowest}$\sc\,$) {\it weight vector} \vs{-4pt}if
$$
(\LIE)_0v\subset\C v, \mbox{ \ \ and \ \ } (\LIE)_+v=0\mbox{ \
(\,resp., }
 (\LIE)_-v=0\,\mbox{)}.
\vs{-4pt}$$%
\par
We divide the proof of Theorem 2.2 into 3 lemmas. \par
{L$\sc\rm EMMA$ 3.1.} \ {\it Suppose $V$ is an irreducible
quasifinite \LIE-module without highest and lowest weight vectors.
Then $t|_{V_i}:V_i\to V_{i+N}$ and $t^{-1}|_{V_i}:V_i\to V_{i-N}$
are injective and thus bijective for all $i\in\Z$ (recall (2.4)
that $t$ has degree $N$). In particular, by letting $K={\rm
max}\{{\rm dim\,}V_p\,|\,p\in[1,N]\}$, we have ${\rm dim\,}V_i\le
K$ for
$i\in\Z$; thus $V$ is uniformly bounded.}%
\par%
{\it Proof.} \ Say $tv_0=0$ for some $0\ne v_0\in V_i$. By
shifting the grading index if necessary, we can suppose $i=0$.
Since $x|_{V_0}:V_0\to V_{[0,N)}$ for $x\in(\LIE)_{[0,N)}$ are
linear maps on finite dimensional vector spaces, there exists a
finite subset $S\subset(\LIE)_{[0,N)}$ such that all
$x|_{V_0},x\in(\LIE)_{[0,N)}$ are linear combinations of
$S|_{V_0}=\{y|_{V_0}\,\,\big|\,\,y\in S\}$. This implies that
$(\LIE)_{[0,N)}v_0=({\rm span\,}S)v_0$. Recall that $\ad_t$ is
locally nilpotent such that
$(\LIE)_j\subset\ad_t^k((\LIE)_{[0,N)})$ for $j>0$, where $k\ge0$
is the integer such that $0\le j-kN<N$. Choose $n>0$ such that
$\ad_t^n(S)=0$. Let $j\ge nN$. Then $k\ge n$ and we \vs{-4pt}have
$$
(\LIE)_jv_0\subset(\ad_t^k((\LIE)_{[0,N)}))v_0=t^k(\LIE)_{[0,N)}v_0=t^k({\rm
span\,}S)v_0 =(\ad_t^k({\rm span\,}S))v_0=0. %
\vs{-4pt}$$%
 This
means that $(\LIE)_{[nN,\infty)}v_0=0.$\par%
The rest of the proof
is exactly like that of Proposition 2.1 in [13]. %
\qed\par%
{L$\sc\rm EMMA$ 3.2.} \ {\it A nontrivial irreducible uniformly
bounded module $V$ is a module $V(\a)$ or $\ol V(\a)$ for some
$\a\in\C$.}\par%
{\it Proof.} \ Let $V'={\rm span}\{v\in V\,|\,\HH v\subset\C v\}$
be the space spanned by weight vectors. Clearly $V'$ is a
submodule. Since ${\rm dim\,}V_i<\infty$, there exists at least a
common eigenvector (i.e., a weight vector) of $\HH$ in $V_i$ for
$i\in\Z$, i.e., $V'\ne0$. Thus $V=V'$ is a weight module. Since
$\cc,I$ are central elements, we \vs{-3pt}have (cf.~Remark 3.4
below)
$$
\cc|_V=c_0{\sc\!}\cdot{\sc\!}{\bf1}_V,\;\;\;\;I|_V=c_1{\sc\!}\cdot{\sc\!}{\bf1}_V,
$$%
\vs{-3pt}for some $c_0,c_1\in\C$. %
Let%
$$ \Vir={\rm span}\{t^iD,\cc\,|\ i\in\Z\},\,\;\;\WW={\rm
span}\{t^iD^j,\cc\,|\ i\in\Z,\,j\in\Z_+\}, %
$$%
be subalgebras of \LIE\ isomorphic to the Virasoro algebra and the
Lie algebra $\wh{\DD^1}$ respectively. For $m\in[0,N-1]$,
\vs{-3pt}set
$$
V[m]=\bigoplus_{i\in\Z}V_{iN+m}. \eqno(3.1)
$$%
Clearly, $V[m]$ is a uniformly bounded module over $\Vir$ or
$\WW$, and $V=\oplus_{m=0}^{N-1}V[m]$. Since $V\ne0$, we have
$V[m]\ne0$ for some $m$. Say, $V[0]\ne0$. Obviously, a
composition factor of the
$\Vir$-module $V[0]$ is a $\Vir$-module of
the intermediate series, on which the central element $\cc$ must
act
trivially (cf.~[12, 14]); thus $c_0=0$ (and so we can omit $\cc$ in the following discussion).%
\par%
By the structure of uniformly bounded $\WW$-modules in [13], we
have $c_1=0,\pm1$. If $c_1=0$, by [13], each $V[m]$ must be a
trivial $\WW$-module; thus $V$ is trivial as a $\WW$-module. But
since $[\WW,\LIE]=\LIE$, we obtain that $V$ is a trivial
$\LIE$-module, contradicting the assumption of the lemma.
Thus $c_1\ne0$. If necessary, by using the isomorphism $\si$ in
(2.8) (which interchanges $V(\a)$ with $\ol V(\a)$), we can always
suppose $c_1\ne-1$ (since $\si(I)=-I$). Thus $c_1=1$. Thus again
by [13], each $\WW$-module $V[m]$ must have the form $A_{p,G}$
defined in [13], i.e., there exist $K_m\ge0$ (might depend on $m$)
and a $K_m\times K_m$ diagonal matrix $G_m$ such that we can
choose a suitable basis
$Y_{kN+m}=(y_{kN+m}^{(1)},...,y_{kN+m}^{(K_m)})$ of $V_{kN+m}$
for $k\in\Z$ \vs{-4pt}satisfying%
 $$
(t^iD^j)Y_{kN+m}=Y_{(i+k)N+m}(k+G_m)^j,
\vs{-4pt} \eqno(3.2)$$%
for $i,k\in\Z,\,j\in\Z_+$, where the right-hand side is the
vector-matrix multiplication by regarding $Y_{kN+m}$ as a row
vector, and where, here and below, when the context is clear, in
an expression involving matrices, we always identify a scalar
(like $k$ in the right-hand side) with the corresponding scalar
matrix of the suitable rank.
\par
Note that for $p,q\in[1,N]$, $tE_{p,q}$ has degree $N+p-q$, thus
$(tE_{p,q})V_{kN+m}\subset V_{(k+1)N+m+p-q}$. Let $k_1,m_1$ be the
integer
such that $m_1\in[0,N-1]$ \vs{-4pt}and %
$$kN+m+p-q=k_1N+m_1.\vs{-4pt}$$
 Then we can \vs{-4pt}write
$$
(tE_{p,q})Y_{kN+m}=Y_{(k_1+1)N+m_1}P_{p,q}^{(k,m)}, %
\vs{-4pt}\eqno(3.3)$$%
where $P_{p,q}^{(k,m)}$ is some $K_{m_1}\times K_m$ matrix.
Applying $[t^i,tE_{p,q}]=0$ and $[D,E_{p,q}]=0$ to $Y_{kN+m}$,
using (3.2) and (3.3), we obtain that $P_{p,q}^{(k,m)}$ does not
depend on $k$ (and denote it by $P_{p,q}^{(m)}$) \vs{-4pt}and
$$G_{m_1}P_{p,q}^{(m)}=P_{p,q}^{(m)}G_m.%
\vs{-4pt}\eqno(3.4)$$%
 Applying
$[t^iD^j,tE_{p,q}]=\sum_{s=1}^j(^j_s)t^{i+1}D^{j-s}E_{p,q}$ to
$Y_{kN+m}$, by induction on $j$, we \vs{-4pt}obtain
$$
(t^iD^jE_{p,q})Y_{kN+m}=Y_{(k_1+i)N+m_1}P_{p,q}^{(m)}(k+G_m)^j.
\vs{-7pt}\eqno(3.5)$$%
\par
Let $K=\sum_{m=0}^{N-1}K_m$ and let $U=\oplus_{m=0}^{N-1}V_m$ be a
subspace of $V$ of dimension $K$. Then
$$%
(Y_0,Y_1,...,Y_{N-1})=
(y_0^{(1)},...,y_0^{(K_0)},y_1^{(1)},...,y_1^{(K_1)},...,
y_{N-1}^{(1)},...,y_{N-1}^{(K_{N-1})}),$$%
 is a basis of $U$. For
each $E_{p,q}\in gl_N$, we define a linear transformation
$\rho(E_{p,q})$ of $U$ as follows:
$$
\rho(E_{p,q})Y_m=Y_{m_1}P_{p,q}^{(m)}, \eqno(3.6)
$$
for $m\in[0,N-1]$, where $m_1\in[0,N-1]$ such that $m_1\equiv
m+p-q\,({\rm mod\,}N)$.
 This uniquely defines a linear map
$\rho:gl_N\to{\rm End}(U)$. We prove that
$\rho(E_{p,q}E_{p',q'})=\rho(E_{p,q})\rho(E_{p',q'}).$ By shifting
the index of $V_k$ if necessary, it suffices to prove
$$\rho(E_{p,q})\rho(E_{p',q'})Y_0
=\rho(E_{p,q}E_{p',q'})Y_0=\d_{q,p'}\rho(E_{p,q'})Y_0.
\eqno(3.7)$$%
Note from (2.1) and (2.2) that%
 $$[DE_{p,q},tE_{p',q'}]=
tE_{p,q}E_{p',q'}+tD[E_{p,q},E_{p',q'}]=\d_{q,p'}tE_{p,q'}+
\d_{q,p'}tDE_{p,q'}-\d_{q',p}tDE_{p',q}.%
\eqno(3.8)$$%
First suppose $p\ge q,\,p'\ge q'$, applying (3.8) to $Y_{kN}$, by
(3.5), we obtain $$%
P_{p,q}^{(p'-q')}(\ol
k+1)P_{p',q'}^{(0)}-P_{p',q'}^{(p-q)}P^{(0)}_{p,q}\ol k
=\d_{q,p'}P_{p,q'}^{(0)}+(\d_{q,p'}P_{p,q'}^{(0)}-\d_{q',p}P_{p',q}^{(0)})\ol
k,%
\eqno(3.9)$$%
where, we denote $\ol k=k+G_0$. Since $\ol k$ commutes with
$P_{p',q'}^{(0)}$ by (3.4), regarding expressions in (3.9) as
polynomials on $\ol k$, by comparing the coefficient of $\ol k^0$,
we obtain
$$P_{p,q}^{(p'-q')}P_{p',q'}^{(0)}=\d_{q,p'}P_{p,q'}^{(0)},%
\eqno(3.10)$$%
which is equivalent to (3.7). By symmetry, we also have (3.7)
if $p\le p,\,p'\le q'$. Now suppose $p<q,\,q'\ge q'$. Again apply
(3.8) to $Y_{kN}$, we have (3.9) with $P_{p',q'}^{(p-q)}$
replaced by $P_{p',q'}^{(N+p-q)}$. Thus we still have (3.10)
which is again equivalent to (3.7). Finally suppose $p\ge
q,\,p'<q'$. Then we have (3.9) and (3.10) with
$P_{p,q}^{(p'-q')}$ replaced by $P_{p,q}^{(N+p'-q')}$ and again we
have (3.7).
\par
 Thus $\rho$ is a representation of the simple associative
algebra $gl_N$ ($\,={\rm End}_N\,$). Thus
$U=\oplus_{s=1}^{n}U^{(s)}$ is decomposed as a direct sum of
simple $gl_N$-submodules $U^{(s)}$ such that each $U^{(s)}$ is
either the natural $gl_N$-module ($\,\cong\C^N\,$) or the trivial
module. Since $I|_U$ is the identity map, we have that each
$U^{(s)}$ is the natural module. Since $[D,gl_N]=0$ and $D|_U$ is
diagonalizable, we can choose submodules $U^{(s)}$ such that
$D(U^{(s)})\subset U^{(s)}$ and $D|_{U^{(s)}}$ is a scalar map.
Now clearly $U^{(1)}$ generates a simple \LIE-submodule of $V$ of
the form $V(\a)$ (cf.~(3.5) and (3.6)). Since $V$ is irreducible,
we have $V=V(\a)$ for some $\a\in\C$ (if we have used the
isomorphism $\si$ in (2.8) in the
above proof, then $V$ is the module $\ol V(\a)$).%
\qed\par%
{L$\sc\rm EMMA$ 3.3.} \ {\it A nontrivial indecomposable uniformly
bounded module $V$ is a module of the form $V(m,\a)$ or $\ol
V(m,\a)$.}%
\par{\it Proof.} \
First note that a central element,
while not necessarily acting by a scalar on an indecomposable module,
nevertheless has only one eigenvalue (cf.~Remark 3.4 below). Let $c_0$
and $c_1$ be the eigenvalue of $\cc$ and $I$ respectively. As in
the arguments of the proof of Lemma 3.2, we have $c_0=0$ (thus
$\cc$ acts nilpotently on $V$) and we can suppose $c_1=1$ (by
making use of the isomorphism $\si$ in (2.8)). Thus each
composition factor of $V$ has the form $V(\a)$. Therefore $V$ has
a finite number $m$ of composition factors ($\,m={\rm
dim\,}V_0\,$). By induction on $m$, it suffices to consider the
case when $m=2$. In this case $V$ is not irreducible but
indecomposable.
\par%
First suppose
$$%
\cc|_V=0,\;\;\;\;I|_V={\bf1}_V.
\eqno(3.11)$$%
Following the proof of Lemma 3.2 (now $G_m$ is not necessarily
diagonal), we have $U=U^{(1)}\oplus U^{(2)}$, and both $U^{(1)}$
and $U^{(2)}$ are the natural $gl_N$-modules. Since $D(U)\subset
U$, $[D,gl_N]=0$ and $V$ is not irreducible, the subspace
$U'=\{u\in U\,|\,Du\in\C u\}$ of eigenvectors of $D$ is a proper
(and thus simple) $gl_N$-submodule of $U$ (isomorphic to $\C^N$ as
a $gl_N$-module) and $D|_{U'}$ is a scalar map $\l$ for some
$\l\in\C$. Thus $U=U'\oplus U''$, where $U''$ is another copy of
$U'$ such that $Du''=\l u''+u'$ for $u''\in U''$, where $u'\in U'$
is the corresponding copy of $u''$. Therefore
$U\cong\C^N\otimes\C^2$ such that $gl_N$ acts on $\C^N$ and $D$
acts on $\C^2$ (and $\a=D|_{\C^2}$ is an indecomposable linear
transformation on $\C^2$), and we obtain $V=V(2,\a)$ (if we have
used the isomorphism $\si$ in (2.8) in the above proof, then $V$
is the module $\ol V(2,\a)$).
\par%
It remains to prove (3.11). We shall use (2.3), which can be
rewritten as follows:%
\vs{-4pt}$$
\begin{array}{ll}
&
[t^{i+j}({d\over dt})^jA,t^{k+l}({d\over dt})^lB]\vs{4pt}\\%
=\!\!& %
{\dis\sum_{s=0}^j}\biggl(\!\begin{array}{c}j\\
s\end{array}\!\biggr)[k+l]_st^{i+j+k+l-s}({d\over dt})^{j+l-s}AB-
{\dis\sum_{s=0}^l}\biggl(\!\begin{array}{c}l\\
s\end{array}\!\biggr)[i+j]_st^{i+j+k+l-s}({d\over dt})^{j+l-s}BA
\vs{4pt}\\
&+\,\d_{i,-k}(-1)^jj!l!
\biggl(\!\begin{array}{c}i+j\\
j+l+1\end{array}\!\biggr){\rm tr}(AB)\cc, \end{array}%
 \vs{-7pt}\eqno(3.12)$$%
where $[k]_j=k(k-1)\cdots(k-j+1)$ is a similar notation to $[D]_j$
in (2.3).%
\par
 Consider the $\WW$-module $V[0]$ (cf.~(3.1)).
As in the proof of Proposition 2.2 of [13], we can choose a basis
$X_0=(x^{(1)}_0,x^{(2)}_0)$ of $V_0$ and define a basis
$X_n=(x^{(1)}_n,x^{(2)}_n)$ of $V_{nN}$ by induction on $|n|$ such
that $tX_n=X_{n+1}$ for $n\in\Z.$
 Assume \vs{-4pt}that
$$\cc X_n=X_nC_n,\,\;\;\;(t^{i+j}(\mbox{${d\over dt}$})^j)X_n=X_{n+i}P_{i,j,n},%
\vs{-4pt}$$%
 for some $2\times2$ matrices $C_n,P_{i,j,n}$.
Using $[\cc,t]=0$, we obtain $C_n=C_0$. Using $[t^{i+j}({d\over
dt})^j,t] =jt^{i+j}({d\over dt})^{j-1}$, we obtain
$P_{i,j,n+1}-P_{i,j,n}=jP_{i+1,j-1,n}$. Thus induction on $j$
\vs{-4pt}gives
$$
\begin{array}{c}
P_{i,0,n}=P_i,\;\;\;P_{i,1,n}=  \bar nP_{i+1}+Q_i,
\vs{4pt}\\
P_{i,2,n}=[\bar n]_2P_{i+2}+2\bar nQ_{i+1}+R_i,\;\;\;
P_{i,3,n}=[\bar n]_3P_{i+3}+3[\bar n]_2Q_{i+2}+3\bar nR_{i+1}+S_i,
\end{array}
\vs{-4pt}$$%
for some $2\times2$ matrices $P_i,Q_i,R_i,S_i$, where $\bar
n=n+G_0$ for some fixed $2\times2$ matrix $G_0$ (cf.~(3.2)),
$[\bar n]_j$ is again a similar notation to $[D]_j$, and $Q_0=0$.
(We use notation $\bar n=n+G_0$ in order to ensure that $Q_0=0$.
Note from $[t{d\over dt},t^{i+j}({d\over dt})^j]=it^{i+j}({d\over
dt})^j$ that $G_0$ commutes with all matrices
$C_0,P_i,Q_i,R_i,S_i$.)%
 \par%
By choosing a composition series of $V[0]$, we can assume that all
these matrices are upper triangular matrices. Furthermore, by the
structure of modules of the intermediate series, we see that %
$$%
P_i\mbox{ have the form }(^{1\ *}_{0\ 1}), \mbox{ and
}C_0,Q_i,R_i,S_i\mbox{ have the form }(^{0\ *}_{0\ 0}).%
\eqno(3.13)$$%
 Thus all
matrices in (3.13) are commutative. %
\par%
Applying $[t^{i+1}{d\over dt},t^k]=kt^{i+k}-\d_{i,-k}(^{i+1}_{\
2})NC$ (cf.~(3.12)) to $X_n$, we obtain
$kP_{i+1}P_k=kP_{i+k}-\d_{i,-k}(^{i+1}_{\ 2})NC_0$, from which we
obtain $P_i=P_0^{1-i}=(1-i)P_0+i$ for $i\in\Z$ (using (3.13), we
have $(P_0-1)^2=0\,$)
 and $C_0=0$ (thus $\cc|_{V[0]}=0$ and similarly $\cc|_{V[m]}=0$
 for $m\in[0,N-1]$,
and so in the following, we can omit $\cc$). Applying
$[t^2({d\over dt})^2,t^k]=2kt^{k+1}{d\over dt}+[k]_2t^k$ to $X_n$,
by comparing the coefficient of $\ol n^0$, we obtain
$([k]_2P_2+2kQ_1)P_k=2kQ_k+[k]_2P_k.$ \vs{-4pt}Thus%
$$Q_k=Q_1P_k+\mbox{${k-1\over2}$}(P_2-1)P_k=Q_1+\mbox{${k-1\over2}$}(1-P_0),
\vs{-4pt}\eqno(3.14)$$%
if $k\ne0$ (using (3.13), we have
$Q_1P_k=Q_1,\,(P_2-1)P_k=P_2-1\,$). Letting $j=2,l=0$ and applying
(3.12) to $X_n$, we see that (3.14)
 also holds for $k=0$. Since $Q_0=0$, we obtain
 $Q_k={k\over2}(1-P_0)$. Similarly, Letting $i=l=0,j=3$ and
 applying (3.12) to $X_n$, we obtain
 $([k]_3P_3+3[k]_2Q_2+3kR_1)P_k=3kR_k+3[k]_2Q_k+[k]_3P_k,$ from
 which we \vs{-4pt}obtain
 $$%
 R_k=R_1+\mbox{${(k-1)(k-2)\over3}$}(P_3-1)+(k-1)(Q_2-Q_k)
 =R_1+\mbox{${(k-1)(k-2)\over6}$}(1-P_0).
 \vs{-4pt}$$%
 Finally letting $j=3,l=0$ and applying (3.12) to $X_n$, we obtain $P_0=1$.
  Thus $I|_{V[0]}={\bf1}_{V[0]}$. Similarly $I|_{V[m]}={\bf1}_{V[m]}$.
 This
 proves (3.11), thus the lemma.
\qed\par%
Theorem 2.2 now follows from Lemmas 3.1-3.
\par%
{R$\sc\rm EMARK$ 3.4.} \ We would like to point out that a central
element does not necessarily act as a scalar on an indecomposable
module since we do not assume that a central element acts
diagonalizably. Thus there is a gap in the assertion in (2.1) of
[13]. This gap has been filled in the above proof of Lemma 3.3.
\vs{10pt}\par%
 \cl{REFERENCES}
\lineskip=3.pt \parskip=0.03truein \small
\par\ni\hi3.5ex\ha1
[1] B.~Bakalov, V.~G.~Kac, A.~A.~Voronov, Cohomology of conformal
algebras, {\it Comm.~Math.~Phys.} {\bf200} (1999), 561-598.
\par\ni\hi3.5ex\ha1
[2] C.~Boyallian, V.~Kac, J.~Liberati and C.~Yan,
    Quasifinite highest weight modules of the Lie algebra of matrix
    differential operators on the circle,
    {\it J.~Math.~Phys.} {\bf39} (1998), 2910-2928.
\par\ni\hi3.5ex\ha1
[3] V.~Chari, Integrable representations of affine Lie algebras,
    {\it Invent.~Math.} {\bf85} (1986), 317-335.
\par\ni\hi3.5ex\ha1
[4] E.~Frenkel, V.~Kac, R.~Radul and W.~Wang,
    $\WW_{1+\infty}$ and $\WW(gl_N)$ with central charge $N$,
    {\it Comm. Math.~Phys.} {\bf170} (1995), 337-357.
\par\ni\hi4ex\ha1
[5] V.~G.~Kac, {\it Vertex algebras for beginners}, American
Mathematical Society, Providence, 1996.
\par\ni\hi3.5ex\ha1
[6] V.~G.~Kac, The idea of locality, in ``Physical applications
and mathematical aspects of geometry, groups and algebras'', H.-D.
Doebner et al, eds., World Sci., Singapore, 1997, 16-32. 
\par\ni\hi3.5ex\ha1
[7] V.~G.~Kac, Formal distribution algebras and conformal
algebras, a talk at the Brisbane, in Proc.~XIIth International
Congress of Mathematical Physics (ICMP '97) (Brisbane), 80-97.
\par\ni\hi3.5ex\ha1
[8] V.~G.~Kac and A.~Radul, Quasi-finite highest weight modules
over the
    Lie algebra of differential operators on the circle,
    {\it Comm.~Math.~Phys.} {\bf157} (1993), 429-457.
\par\ni\hi3.5ex\ha1
[9] V.~G.~Kac and A.~Radul, Representation theory of the vertex
algebra
     $\WW_{1+\infty}$, {\it Trans. Groups} {\bf1} (1996), 41-70.
\par\ni\hi4ex\ha1
[10] V.~G.~Kac, W.~Wang and C.~H.~Yan, Quasifinite representations
of classical
     Lie subalgebras of $\WW_{1+\infty}$, {\it Adv.~Math.}
     {\bf139} (1998), 46-140.
\par\ni\hi4ex\ha1
[11] W.~Li, 2-Cocycles on the algebra of differential operators,
     {\it J.~Alg.} {\bf122} (1989), 64-80.
\par\ni\hi4ex\ha1
[12] O.~Mathieu, Classification of Harish-Chandra modules over the
     Virasoro Lie algebra, {\it Invent.~Math.} {\bf 107} (1992), 225-234.
\par\ni\hi4ex\ha1
[13] Y.~Su, Classification of quasifinite modules over
        the Lie algebras of Weyl type, {\it Adv.~Math.} {\bf174} (2003),
        57-68.
\par\ni\hi4ex\ha1
[14] Y.~Su, Classification of Harish-Chandra modules over the
     higher rank Virasoro algebras, {\it Comm. Math.~Phys.}
{\bf 240} (2003), 539-551.
\par\ni\hi4ex\ha1
[15] Y.~Su and K.~Zhao, Isomorphism classes and automorphism
groups of
     algebras of Weyl type, {\it Science in China A} {\bf45}
     (2002), 953-963.
\par\ni\hi4ex\ha1
[16] X.~Xu, Equivalence of conformal superalgebras to Hamiltonian
superoperators, {\it Alg.~Colloq.} {\bf8} (2001), 63-92.
\par\ni\hi4ex\ha1
[17] X.~Xu, Simple conformal algebras generated by Jordan
algebras, preprint, math.QA/0008224. 
\par\ni\hi4ex\ha1
[18] X.~Xu, Simple conformal superalgebras of finite growth, {\it
Alg.~Colloq.} {\bf7} (2000), 205-240.
\par\ni\hi4ex\ha1
[19] X.~Xu, Quadratic Conformal Superalgebras, {\it J.~Alg.}
{\bf231} (2000), 1-38.
\end{document}